\documentclass[twoside,10pt]{article}
\usepackage{amssymb,amsmath,euscript}
\usepackage{graphicx,color,caption}

\usepackage{algorithm,algpseudocode}

\newtheorem{proposition}{Proposition}[section]
\newtheorem{theorem}[proposition]{Theorem}

\newtheorem{definition}[proposition]{Definition}

\newcommand{\qed}{\hphantom{.}\hfill $\Box$\medbreak}

\def\S{\mathbb{S}}
\def\P{\mathcal{P}}
\def\D{\mathcal{D}}

\def\T{\mathbb{T}}
\def\Q{\mathcal{Q}}
\def\O{\mathcal{O}}
\def\A{\mathcal{A}}
\def\J{\mathcal{J}}
\def\B{\mathcal{B}}
\def\C{\mathcal{C}}

\def\I{\mathcal{I}}
\def\J{\mathcal{J}}
\def\x{{\bf x}}
\def\u{{\bf u}}
\def\0{{\bf 0}}
\title{\bf{Characterization Tensors of Balanced Incomplete Block Designs}
\thanks{This research was supported by the National Natural Science Foundation of China (11301022,11431002), the State Key Laboratory of Rail Traffic Control and Safety, Beijing Jiaotong University (RCS2014ZT20, RCS2014ZZ01), and the Hong Kong Research Grant Council (Grant No.
PolyU 502111, 501212, 501913 and 15302114).}}
\author{Liqun Qi \thanks{Department of Applied Mathematics, The Hong Kong Polytechnic University, Hung Hom, Kowloon, Hong Kong. ({\tt
liqun.qi@polyu.edu.hk}).}  \hspace{1mm}   Ziyan
Luo,  \thanks{ State Key Laboratory of Rail Traffic Control and Safety, Beijing
Jiaotong University, Beijing 100044, P.R. China; ({\tt
starkeynature@hotmail.com}).}
}

\begin{document}
\maketitle

\begin{abstract}
Balanced incomplete block designs (BIBDs) have wide applications in engineering, business and sciences.   In this paper, for each $(v, k, \lambda)$-BIBD, we construct  a strongly symmetric $k$-th order $v$-dimensional tensor.    We call such a strongly symmetric tensor the characterization tensor of that BIBD, and the absolute value tensor of
the characterization tensor the signless characterization tensor of that BIBD.   We study some spectral properties of such characterization tensors and signless characterization tensors.  In this way, we provide a new tool to study BIBDs.


\end{abstract}


\section{Introduction}
Balanced incomplete block designs (BIBDs) have wide applications in engineering, business and sciences \cite{Wa2007, Wi2015}.

Given a finite set $X$ of points and integers $k, r, \lambda \ge 1$, we may define a {\bf balanced incomplete block design} $B$ to be a family of $k$-element subsets of $X$, called {\bf blocks} such that for each $i \in X$, there are exactly $r$ blocks containing $i$, and for each pair $(i, j)$ of points in $X$, there are exactly $\lambda$ blocks containing $i$ and $j$.    Let the cardinality of $X$ to be $v$ and the cardinality of $B$ be $b$.   Then we have
$$bk = vr$$
and
$$\lambda(v-1) = r(k-1).$$
We assume the design is simple, i.e., repeated blocks are not allowed.   We also assume that $v > k$.  In design, this is the meaning of the word ``incomplete'' \cite{Wa2007, Wi2015}.

In this paper, for each $(v, k, \lambda)$-BIBD, we construct a strongly symmetric $k$-th order $v$-dimensional tensor.    We call such a strongly symmetric tensor the characterization tensor of that BIBD, and the absolute value tensor of
the characterization tensor the signless characterization tensor of that BIBD.   We study some spectral properties of such characterization tensors and signless characterization tensors.  In this way, we provide a new tool to study BIBDs.

Some notations that will be used throughout the paper are listed here. The $v$-dimensional real Euclidean space is denoted by $\Re^v$, where $v$ is a given natural number. The nonnegative orthant in $\Re^v$ is denoted by $\Re^v_+$, with the interior $\Re^v_{++}$ consisting of all positive vectors. 
Denote $[v]:=\{1,2,\ldots, v\}$. Vectors are denoted by bold letters such as $\x$, $\u$, matrices are denoted by capital letters such as $A$, $P$, and tensors are written as calligraphic capital letters such as $\A$, $\B$. The space of all real $k$th order $v$-dimensional tensors is denoted by $\T_{k,v}$, and the space of all symmetric tensors in $\T_{k,v}$ is denoted by $\S_{k,v}$. 
For a subset $\Gamma\subseteq [v]$, $|\Gamma|$ stands for its cardinality.



\section{Strongly Symmetric Tensors and Eigenvalues}
Let $\A=\left(a_{i_1\ldots i_k}\right)$ be a $k$th order $v$-dimensional real tensor. $\A$ is called a \emph{symmetric} tensor if the entries $a_{i_1\ldots i_k}$ are invariant under any permutation of their indices for all $i_j\in [v]$ and $j\in [k]$, denoted as $\A\in \S_{k,v}$. A symmetric tensor $\A$ is said to be positive semidefinite (definite) if $\A \x^k:=\sum\limits_{i_1,\ldots,i_k\in [v]}a_{i_1\ldots i_k}x_{i_1}\cdots x_{i_k}\geq 0 (>0)$ for any $\x\in \Re^v\setminus \{\0\}$ \cite{Qi2005}. Here, $\x^k$ is a rank-one tensor in $\S_{k,v}$ defined as $\left(\x^k\right)_{i_1\ldots i_k}:=x_{i_1}\cdots x_{i_k}$ for all $i_1$, $\ldots$, $i_k\in [v]$. Evidently, when $k$ is odd, $\A$ could not be positive definite and $\A$ is positive semidefinite if and only if $\A=\O$, where $\O$ stands for the zero tensor. A tensor $\A\in \S_{k,v}$ is said to be (strictly) copositive if $\A x^k\geq 0$ ($>0$) for all $x\in\Re^n_{+}\setminus \{\0\}$ \cite{Qi2013}. The definitions on eigenvalues of symmetric tensors are as follows.

\begin{definition}[\cite{Qi2005}] Let $\A\in \S_{k,v}$ and $\mathbb{C}$ be the complex field. We say that $(\mu, \x)\in {\mathbb{C}}\times \left({\mathbb{C}}^v\setminus \{\0\}\right)$ is an eigenvalue-eigenvector pair of $\A$ if
$\A \x^{k-1}=\mu x^{[k-1]},$ where $\A \x^{k-1}$ and $\x^{[k-1]}$ are all $v$-dimensional column vectors given by
$$\left(\A \x^{k-1}\right)_i:=\sum\limits_{i_2,\ldots,i_k\in[v]} a_{ii_2\ldots i_k} x_{i_2}\cdots x_{i_k},~\left(\x^{[k-1]}\right)_i=x_i^{k-1},~~\forall i\in [v].\eqno(2.1)$$  If the eigenvalue $\mu$ and the eigenvector $\x$ are real, then $\mu$ is called an $H$-eigenvalue of $\A$ and $\x$ an $H$-eigenvector of $\A$ associated with $\mu$. If $\x\in \Re^n_+ (\Re^n_{++})$, then $\mu$ is called an $H^+ (H^{++})$-eigenvalue of $\A$.  The maximum modulus of the eigenvalues of $\A$ is called the spectral radius of $\A$ and denoted by $\rho(\A)$.   The largest $H$-eigenvalue of $\A$ is denoted as $\mu_{\max}(\A)$.  The set of all the $H$-eigenvalues of $\A$ is called the $H$-spectrum of $\A$.
\end{definition}

When $k$ is even, $\A\in \S_{k,v}$ always has $H$-eigenvalues \cite{Qi2005}.   A nonnegative tensor always has $H$-eigenvalues \cite{CQZ2013}.

\begin{definition}[\cite{Qi2005}] Let $\A\in \S_{k,v}$ and $\mathbb{C}$ be the complex field. We say that $(\mu,\x)\in {\mathbb{C}}\times \left({\mathbb{C}}^v\setminus \{\0\}\right)$ is an $E$ eigenvalue-eigenvector pair of $\A$ if
$\A \x^{k-1}=\mu \x$ and $\x^{T}\x=1$, where $\A \x^{k-1}$ is defined as in (2.1). If the $E$-eigenvalue $\mu$ and the eigenvector $\x$ are real, then $\mu$ is called a $Z$-eigenvalue of $\A$ and $\x$ a $Z$-eigenvector of $\A$ associated with $\mu$.
\end{definition}

$\A\in \S_{k,v}$ always has $Z$-eigenvalues \cite{Qi2005}.

A $k$th order $v$-dimensional real tensor $\A=\left(a_{i_1\ldots i_k}\right)$ is called a \emph{strongly symmetric} tensor if $a_{i_1\ldots i_k} = a_{j_1\ldots j_k}$ as long as
$\{ i_1, \ldots, i_k \} = \{ j_1, \ldots, j_k \}$ for all $i_l, j_l\in [v]$ and $l\in [k]$.   For example, for a fourth order $v$-dimensional real tensor $\A=\left(a_{i_1i_2i_3i_4}\right)$, we always have $a_{ijjj}=a_{ijij}=a_{jiii}$ for all $i, j \in [v]$.  Denote the set of all $k$th order $v$-dimensional real strongly symmetric tensors as $\S\S_{k,v}$.  Strongly symmetric tensors were introduced in \cite{QXX2014} and further studied in \cite{QWC2015}.

Let $\A=\left(a_{i_1\ldots i_k}\right) \in \T_{k,v}$.    We call $a_{i_1\ldots i_k}$ a diagonal entry of $\A$ if $i_1= \ldots = i_k$.   Otherwise, we call it an off-diagonal entry of $\A$.
If $| \{ i_1, \ldots, i_k \} | = 2$, we say that $a_{i_1\ldots i_k}$ is a sub-diagonal entry of $\A$.   A tensor in $\T_{k,v}$ is called a diagonal tensor if all of its off-diagonal entries are $0$.  A diagonal tensor in $\T_{k,v}$ is called an identity tensor of $\T_{k,v}$, and denoted by $\I$ if all of its diagonal entries are $1$.  A tensor in $\T_{k,v}$ is called a sub-diagonal tensor if all of its entries, which are not sub-diagonal entries are $0$.   A sub-diagonal tensor in $\T_{k,v}$ is called an identity sub-diagonal tensor of $\T_{k,v}$, and denoted by $\J$ if all of its sub-diagonal entries are $1$ and the other entries are $0$.  Clearly, diagonal tensors and sub-diagonal tensors are strongly symmetric tensors.

\section{Characterization Tensors of Balanced Incomplete Block Designs}

We may regard $(v, k, \lambda)$-BIBD as a $k$-uniform hypergraph $G = (X, B)$ \cite{CD2012, Qi2014}, where $X$ is the vertex set and $B$ is the edge set.  For hypergraphs, for each vertex $i \in X$, let $X_i = \{ e \in B : i \in e \}$.  Then $d(i) = |X_i|$ is called the degree of $i$.   For each pair of vertices $i, j \in X$, let $X_{i, j} = \{ e \in B : i, j \in e \}$.  Then $d(i, j) = |X_{i, j}|$ is called the co-degree of $(i, j)$.   A $k$-uniform hypergraph is called a $r$-regular $k$-uniform hypergraph if the degrees of its vertices are the same as $r$ \cite{CD2012, QSW2014}.  Thus, a $(v, k, \lambda)$-BIBD is corresponding to a $r$-regular $k$-uniform hypergraph with the co-degrees of its vertices are the same as $\lambda$.

For a $k$-uniform hypergraph $G = (X, B)$, there are several tensors associated with it.  The adjacency tensor $\A$ of $G$ is defined as a tensor $\A=\left(a_{i_1\ldots i_k}\right) \in \S_{k,v}$, with $a_{i_1\ldots i_k} = {1 \over (k-1)!}$ if $(i_1, \ldots, i_k) \in B$, and $0$ otherwise \cite{CD2012}.   The degree tensor $\D$ of $G$ is defined as a diagonal tensor in $\S_{k,v}$ with its $i$th diagonal entries as $d(i)$.   The Laplacian tensor and signless Laplacian tensor of $G$ are defined as $\D - \A$ and $\D + \A$ respectively \cite{Qi2014}.  Spectral hypergraph theory via tensors deal with the spectral properties of adjacency tensors, Laplacian tensors and signless Laplacian tensors, and their relations with hypergraphs \cite{CCL2015, CD2012, HQX2015, Qi2014, QSW2014, SSW2015}.  We now define the co-degree tensor $\C=(c_{i_1\ldots i_k})$ of $G$ by $c_{i_1\ldots i_k}= {d(i, j) \over 2^{k-1}-1}$ if $\{ i_1, \ldots, i_k \} = \{ i, j \}, i \not = j$, and $0$ otherwise.  Thus, $G$ is corresponding to a $(v, k, \lambda)$-BIBD if its degree tensor $\D = r\I$ and its co-degree tensor $\C = {\lambda \over 2^{k-1}-1}\J$.
Clearly, all these tensors involved are strongly symmetric tensors.

Now, for a $(v, k, \lambda)$-BIBD, we may define its {\bf characterization tensor} as
$$\P = 2r(k-1)\I + \C - (k-1)\A$$
and its {\bf signless characterization tensor} as
$$\Q = 2r(k-1)\I + \C  + (k-1)\A.$$

Since $\lambda \ge 1$, the hypergraph $G$ corresponding to a $(v, k, \lambda)$-BIBD is connected.  Recall that \cite{HQX2015} a $k$-uniform hypergraph $G = (X, B)$ is called odd-bipartite if $k$ is even and there is a subset $Y$ of $X$ such that for every $e \in B$, $| e \cap Y|$ is odd.   A $(v, k, \lambda)$-BIBD is called odd-bipartite if $k$ is even and the hypergraph $G$ corresponding to it is odd-bipartite.

We now have the following theorem on the spectral properties of $\P$ and $\Q$ of a $(v, k, \lambda)$-BIBD.

\medskip

\begin{theorem}
Suppose that $k, r, \lambda \ge 1$ and $v \ge k$.   Let $\P$ and $\Q$ be the characterization and signless characterization tensors of a $(v, k, \lambda)$-BIBD.   Then

(a) Any eigenvalue $\mu$ of $\P$ and $\Q$ satisfies
$$| \mu -  2r(k-1) | \le 2r(k-1).$$

(b) $\rho(\Q) = 4r(k-1)$.

(c) $\P$ and $\Q$ are co-positive.   When $k$ is even, they are positive semi-definite.

(d) $\mu_{\max}(\P) \le mu_{\max}(\Q)$.  The equality holds here if and only if $k$ is even and the BIBD is odd-bipartite.   If $k$ is even and the BIBD is odd-bipartite, then $\P$ has $H$-eigenvalues and the $H$-spectrum of $\P$ is equal to the $H$-spectrum of $\Q$.
\end{theorem}

{\bf Proof. }  By the Gershgorin theorem for tensor eigenvalues and the definitions of $\P$ and $\Q$, we have (a).

Let $\u$ be the all $1$ vector in $\Re^v$.  Then we have $\Q \u^{k-1} = 4r(k-1)\u^{[k-1]}$.   By (a) and the theory of nonnegative tensors \cite{CQZ2013}, we have (b).

By (a) and \cite{Qi2013, Qi2005}, we have (c).

With a way similar to \cite{HQX2015, SSW2015}, we may prove (d).
\qed

The question is, can we derive some performance properties from the spectral properties of $\P$ and $\Q$.

\section*{Acknowledgements}
The first author is thankful to Prof. Sanming Zhou, who introduced the topic of block design to him, and to Prof. Linyuan Lu for the discussion on co-degrees.


\end{document}